\documentclass[12pt]{article}

\usepackage{amssymb}
\usepackage{epsfig}
\usepackage{amsmath}
\usepackage{amsthm}
\usepackage{subfigure}
\DeclareMathOperator{\transverse}{\cap\kern-7.75pt\top}
\newtheorem{con}{Conjecture}
\newcommand{\Rank}{\operatorname{rank}}
\newtheorem{thm}{Theorem}

\begin{document}

\title{ Reconstruction from projections using dynamics: Non-Stochastic Case
{\bf }} 
\author{Kevin R. Vixie and Gary L. Sandine} 
\date{Los Alamos National Laboratory}

\maketitle

\begin{abstract}

The problem of determining 3D density fields from single 2D
projections is hopelessly under-determined without additional
assumptions. While parameterized inversions are typically used to
solve this problem, we present theoretical results along a different
route to the elimination of indeterminacy. More specifically, we
consider the case in which radiography is being used to study  objects
or processes evolving over some time interval. This evolution can be
measured at several points, generating a sequence of radiographs. A
priori knowledge of constraints on possible evolution permits us to
rule out sequences which are not consistent with those constraints. If
we now consider the space of possible measurement sequences to be the
data space, we find that, under the influence of the constraints, our
data space has increased in dimension while the dimension of the space
of unknowns has remained the same. When enough measurements have been
made, inversion of a dynamically constrained sequence of single angle
radiographs becomes possible.

\end{abstract}

\newpage \small{\tableofcontents}

\newpage

\section[Introduction \textnormal{\\[1ex] \footnotesize We  outline
the paper and briefly introduce the ideas motivating the results we
obtained.}]{Introduction.}
\label{sec:intro}

Reconstructing a series of 3 dimensional density distributions from a
finite number of 2 dimensional measurements is impossible unless prior
assumptions of some sort are used ~\cite{Smith-keenan}. The difficulty
comes from the fact that, without fairly strict assumptions, many
different density fields project to the same radiograph. To state this
another way, a radiographic measurement device, thought of as a
projection operator, has a nontrivial null space. Our approach in this
paper is to discretize the object space and the radiograph
(measurement) space. We then combine a sequence of radiographic
measurements into one super-measurement. Combined with the operator
which determines dynamics, the single time projection operator can be
turned into an extended projection operator that maps a sequence of
objects into a super-measurement. Due to the dynamical constraints,
the dimension of the object sequence space does not grow as the length
of that sequence increases. On the other hand, the size of the data
space (the space of super-measurements) does, implying that
eventually, the extended projection operator has a trivial null space.

\bigskip

Now a look ahead. In section \ref{sec:problem}, we briefly outline the
problem. Section \ref{sec:nota} outlines the notation used for the
rest of the paper. In section \ref{sec:lin} we introduce the problem
in it's linear setting. The key notion turns out to be that of
\emph{transversality} which we recall in section \ref{sec:trans}.  In
this section, we also introduce and prove a theorem which bounds how
slowly the dimension of the projection operator null space decreases
as the number of measurements incorporated in the super-measurement
goes up. There is a  lower bound on the number of measurements that
are needed to get a unique inversion. If $q = $ (dimension of the
object space) $\div$ (dimension of the measurement space), then the
lower bound is simply $\lceil q \rceil$ -- the smallest integer
greater or equal to $q$. We will say that a particular combination of
linear system, $L$ and measurement projection, $P$ has the optimal
reduction property if the number of measurements needed to get an
unique inversion equals this lower bound. Section \ref{sec:opt} looks
at the prevalence of linear systems which (w.r.t. a fixed $P$) having
the optimal reduction property.  Next, in section \ref{sec:nonlin}, we
outline a proof of the extension of one of the results to the case of
nonlinear dynamics. The relation to known results is discussed in
section \ref{sec:RKR}.  We look at a few numerical examples in section
\ref{sec:num} and begin to explore the relationship between
over-determination and noise reduction. We close with a summary and
discussion.

\section[The Problem \textnormal{\\[1ex]\footnotesize Our particular
problem, that of deducing 3-dim density fields from 2-dim radiographic
projections of that field, is introduced here. The abstract
formulation is then given. We arrive at a mathematical problem arising
in many data driven analysis situations. Predictably, it goes by
different names in different contexts.}]{The Problem}
\label{sec:problem}

\begin{figure}[h]
\centerline{ \epsfxsize 5in \epsffile{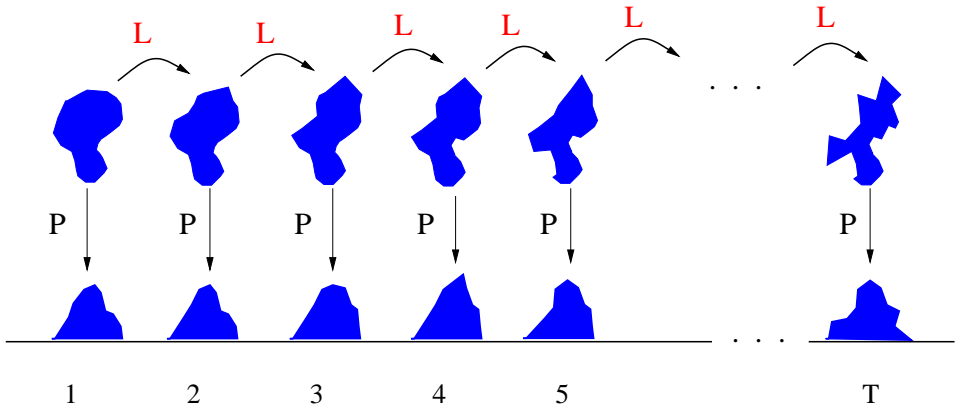}}
\caption{The Problem}
\label{fig:the_problem}
\end{figure}

Radiographic experiments measuring very fast events typically produce
data consisting of a sequence of 2-d projections. These 2-d
projections are created by bombarding some 3-d distribution of density
-- the object -- with penetrating radiation of some sort such as
high-energy x-rays or protons. The number of angles at which the data
is taken is typically 1. We idealize this to get the following model.
The object will be a point $x$ in an object space $X$ which changes
from measurement to measurement as dictated by a linear operator $L$,
acting on $X$. The measurements $d$ which lie in the measurement space
$D$ will be generated by the action of a measurement or projection
operator $P$. Thus, if the object and measurement at time $t \in
\Bbb{N}$ are denoted $x_t$ and $d_t$ respectively, we can express the
actions of the operators $L$ and $P$ in the following way:  $x_{t+1} =
L(x_t)$ and $d_t = P(x_t)$. See figure \ref{fig:the_problem}.

We define the extended (or experimental) spaces to be the product
spaces $\tilde{X} \equiv X^T$ and $\tilde{D} \equiv D^T$ where T is
the number of observation times in a particular experiment.  If we
have a particular sequence of points in the object space, then this
sequence is a single point $\tilde{x} = (x_1,x_2,...,x_T)$ in
$\tilde{X}$. The measurement process produces a point  $\tilde{d}=
(d_1,d_2,...,d_T)$  in the extended (or super-)measurement space
$\tilde{D}$. This can be succinctly expressed using the extended
projection operator $\tilde{P} \equiv P^T$ since then  $\tilde{d}_t =
\tilde{P}(\tilde{x}_t)$.

\bigskip

If A is defined to be the T-1 by T block matrix,

  \[ \left[ \begin{array}{rrrrrr}
                  L & -I &  0 & 0  & ... & 0 \\
                  0 &  L & -I & 0  & ... & 0 \\
                  0 &  0 &  L & 0 & ... & 0 \\
                  ...&...&... &...& ... & ...\\
                  0 &  0 &  0 &   & L   & -I \\
             \end{array} \right] \]

then the null space of this operator, denoted $N_A$, is exactly the
set of elements of $\tilde{X}$ which satisfy the dynamics. That is, $N_A =
\{\tilde{x} \in \tilde{X} | x_{t+1} = L{x_t}\text{ where } \tilde{x} =
(x_1,x_2,...,x_T) \}$.

\bigskip

Let $N_{\tilde{P}}$ be denote the null space of $\tilde{P}$.  The
inverse problem is now solvable when $N_A \cap N_{\tilde{P}} = \{0\}$.

\section[Notation \textnormal{\\[1ex] \footnotesize Oh, the pleasures
and miseries of the use and abuse of creative naming and labeling! We
try to generate a notation that helps the readers in their quest to
follow the derivations. We resist the temptation to make our notation
look good by pointing out examples of notorious notation}]{Notation}
\label{sec:nota}

We now establish the notation we will use throughout, except in
section ~\ref{sec:nonlin} where we find it more convenient to modify
the notation. This section should be used as a reference.

\begin{description}
\item[\hspace{1in}$T \equiv$] The number of radiographs.
\item[\hspace{1in}$X \equiv$] The Space of objects - we will use
$\Bbb{R}^n$.
\item[\hspace{1in}$x \equiv$] An element of $X$.

\item[\hspace{1in}$\tilde{X}\equiv$] $X^T = X \times X\times
...\times X $.

\item[\hspace{1in}$\tilde{x} \equiv$] An element $(x_1,...,x_T)$ of
$\tilde{X}$.

\item[\hspace{1in}$D\equiv$] The space of radiographs - we will use
$\Bbb{R}^m$.

\item[\hspace{1in}$ d \equiv$] An element of $D$.

\item[\hspace{1in}$\tilde{D} \equiv $]$D^T = D\times D\times
...\times D $.

\item[\hspace{1in}$\tilde{d} \equiv$] An element $(d_1,...,d_T)$ of
$\tilde{D}$.

\item[\hspace{1in}$L \equiv$] The linear operator on $X$ that gives
the dynamics: $x_{i+1} = L(x_i)$.

\item[\hspace{1in}$P \equiv$] The projection (measurement) operator
$P:X \rightarrow D$. We assume that $P$ is full rank since otherwise
we may choose a smaller $D$ and consider $P$ with this restricted
range to get a full rank $P$.

\item[\hspace{1in}$N \equiv$] The null space of $P$.

\item[\hspace{1in}$p \equiv$] The dimension of the null space of $P$.
\item[\hspace{1in}$\tilde{P} \equiv$] The extended or product
projection operator $\tilde{P}:\tilde{X} \rightarrow \tilde{D}$.

  \[ \tilde{P} = \left[ \begin{array}{rrrrrr}
                  P & 0  &  0 & 0  & ... & 0 \\
                  0 &  P &  0 & 0  & ... & 0 \\
                  0 &  0 &  P & 0 & ... & 0 \\
                  ...&...&... &...& ... & ...\\
                  0 &  0 &  0 &   & 0   & P \\
             \end{array} \right] \]

\item[\hspace{1in}$A \equiv $] An operator from $\tilde{X}$ to $X^{T-1}$
    
  \[ \left[ \begin{array}{rrrrrr}
                  L & -I &  0 & 0  & ... & 0 \\
                  0 &  L & -I & 0  & ... & 0 \\
                  0 &  0 &  L & 0 & ... & 0 \\
                  ...&...&... &...& ... & ...\\
                  0 &  0 &  0 &   & L   & -I \\
             \end{array} \right] \]

\item[\hspace{1in}$N_A \equiv$] null space of $A$ 
\item[\hspace{1in}$\tilde{P}_{N_A} \equiv$] $\tilde{P}(N_A)$.
\item[\hspace{1in}$ {[}A,...,C{]} \equiv$] The Cartesian product $ A \times ... \times C$.
\end{description}

And we use the following standard notation.

\begin{description}
\item[\hspace{1in}$\text{dim}(H)\equiv$] The dimension of the space or subspace H
\item[\hspace{1in}$M^{\perp} \equiv$] The orthogonal complement of $M$.
\end{description}

And the almost standard notation ...

\begin{description}
\item[\hspace{1in}$S \; \transverse \; Q \equiv $] $S$ intersects $Q$
transversely.
\end{description}

with the slight twist that when $S$ and $Q$ are subspaces, the
intersection always includes the zero vector. So for example, a
transverse intersection of two 1-dimensional subspaces of $R^3$ is
just \{$0$\}, instead of empty intersection as would be the case for
two 1-dimensional curves in $R^3$.

\section[Linear Stationarity implies easy 
solutions \textnormal{\\[1ex] \footnotesize That is, {\it conceptually
easy} ... easy  might be very difficult indeed if the dimensions
involved are big enough! In this section we reduce our problem to
controlling the dimension of a certain sequence of intersections
(which we want to force to 0!) }]{Linear Stationarity implies easy
solutions}
\label{sec:lin}

Suppose now that $X = \Bbb{R}^n$.  Pick a basis for $X$, $\{b_i\}$
where $ i = 1 ... n$. Then $N_A$ is spanned by $\tilde{b}_i \in
\tilde{X}$ where $\tilde{b}_i = (b_i,L(b_i),...,L^{T-1}(b_i)) $. For
example,  suppose that the linear operator $L$ has a complete set of
eigenvectors $\{\omega_i\}$, where $i$ goes from $1$ to $n$. Then
$N_A$ is spanned by $\tilde{\omega}_i \in \tilde{X}$ where
$\tilde{\omega}_i = (\omega_i,L(\omega_i),...,L^{T-1}(\omega_i))
=(\omega_i,\lambda_i\omega_i,...,\lambda_i^{T-1}\omega_i) $.

\bigskip

Now suppose that $D = \Bbb{R}^m$. Form $E_{\tilde{b}}$, the  $n$ by $mT$ matrix
where row $i$ is $\tilde{b}_i = (b_i,L(b_i),...,L^{T-1}(b_i)) $. We
then have that the inverse problem of getting $\tilde{x}$ from
$\tilde{P}(\tilde{x})$ has a unique solution iff the rank of
$E_{\tilde{b}}$ is $n$. (This is virtually identical to the usual test
for observability from control theory, except that I am not assuming
that $T=n$. See for example ~\cite{polderman-willems} p. 178 or
~\cite{sontag-1} p. 271).  We would like to know how big T needs to
be: how many measurements do we need? In the next section we begin to
answer this question with an upper bound on T provided the null space
of $P$ and $L$ do not have a special relation. 

\bigskip

Before we do this, let us look at  $N_A \cap N_{\tilde{P}}$ a little
more carefully.  If the ``true'' sequence in $\tilde{X}$ is
$\tilde{x}^*$ and we have measured $\tilde{d}^* =
\tilde{P}(\tilde{x}^*)$, then anything in $N^T$ can be added to
$\tilde{x}^*$ without changing the observed data
$\tilde{P}(\tilde{x}^*)$. But since we are only interested in those
null sequences that also satisfy the dynamics, we can restrict
ourselves to $n^*$ satisfying,

\[ n^* = (n_1,n_2,...,n_T) \in [N,N \cap L(N), N \cap L(N \cap L(N)), N \cap L( N \cap L( N \cap L(N) )), ... ], \]

where we are using $[A,B,...,C]$ to denote $A\times B \times
... \times C $.  For ease of reference let

\[
\begin{array}{ll}
    N_1 &\equiv N \\
    N_2 &\equiv N \cap L(N_1) \\
    N_3 &\equiv N \cap L(N_2) \\
     ...&    ...\\
    N_T &\equiv N \cap L(N_{T-1}).\\
\end{array}
\]

 This permits us to rewrite the above expression for
$n^*$:

 \[ n^* = (n_1,n_2,...,n_T) \in [N_1,N_2,...,N_T]. \]

\bigskip

Now assume that $L$ is invertible. If $N_T = \{0 \}$, it follows
that $E_{\tilde{b}}$ has rank $= n$ and $N_A \cap N_{\tilde{P}} =
\{0\}$. If $L$ is not invertible, then $N_T = \{0\}$ does not imply
that $N_A \cap N_{\tilde{P}}$ is trivial since $N_A \cap
N_{\tilde{P}}= \{ \text{all possible values of} (L^{-1}(L^{-1}(
... L^{-1}(0)...)),...,L^{-1}(0),0)\}$. Therefore, trivial $N_T$ is a
necessary but insufficient condition for the invertibility of the data
$\tilde{d}^*$ when $L$ is not invertible.

\bigskip

{\bf Remark:} Note that $ n^* = (n_1,n_2,...,n_T) \in
[N,N \cap L(N), N \cap L(N \cap L(N)), ... ] $ is actually more than 
the set of ``null'' solutions. It turns out to be small enough to meet our needs.
 
\section[Transversality is (more than) enough \textnormal{\\[1ex] \footnotesize What conditions on those intersections, which we met in the previous section, guarantee that the dimension of those intersections goes to 0? Transversality does! In fact it implies that the convergence to 0 is as fast as possible! We include a result about how slow the dimension can go to zero and still get there.}]{Transversality is (more than) enough}
\label{sec:trans}

We now study how the dimension of $N_k$ depends on $k$. What
conditions imply that eventually dim($N_k$) becomes $0$ for some $k$?
How prevalent are the $L$'s for a fixed $P$ that have $T^F \equiv$ \{ minimal T such that $N_{T-1} = \{0\}$ \} = $\lceil
n/d \rceil$. That is, how prevalent are the $L$'s having the optimal
reduction property? 

\bigskip

Let us begin by reiterating the
definitions of the $N_k$.

\[
\begin{array}{ll}
    N_1 &\equiv N \\
    N_2 &\equiv N \cap L(N_1) \\
    N_3 &\equiv N \cap L(N_2) \\
     ...&    ...\\
    N_T &\equiv N \cap L(N_{T-1}).\\
\end{array}
\]

If $G$ and $H$ are linear subspaces of $J$ then the only situation
stable to small perturbations is that the intersection of $G$ and $H$
has minimal dimension. That is, we expect dim($G \cap H$) = dim($G$) +
dim($H$) - dim($J$) where a non-positive result indicate the trivial
intersection of dimension zero. If this is the case we shall say that
$G$ and $H$ are transverse and will write this as $G \; \transverse \;
H$.

\bigskip

Referring to the above definitions of the $N_k$'s we see that these sets
decrease in size at the maximum allowable rate if the intersections
defining them are transverse. For example, if as above,dim($X$) = $n$
and dim($N$) = $p$, transverse intersections imply that the sequence
of dimensions is $p$ , $p+p-n$ , $p+p+p-n-n$ ,... or if we note that $d=n-p$
then the sequence is $p$ , $p-d$ , $p-2d$ , $p-3d$ ,... and we get that
dim($N_{\lceil n/d \rceil}$) = $0$ (Remember that we are assuming that $P$ is full rank.)

\bigskip

Suppose that the intersections are not transverse. Then we still have
the following lower bound on the rate at which the dimension of the
$N_i$'s decrease.

\begin{thm}[Minimal Reduction Theorem]
If there is no nontrivial subspace $G$ of $N$ such that $L(G) \subset G$
then for $i \leq \dim N + 1$, dim($N_i$) $\leq$ dim($N$)$-i+1$.
\end{thm} 

\begin{proof}
Assuming for the moment, that $N_{i+1}\subseteq N_i$, we get that $\dim
N_{i+1}\leq\dim N_i$.  Then, if $\dim N_{i+1}=\dim N_i$, we conclude that
$N_i=N_{i+1}$.  By definition of $N_{i+1}$, we then have $N_i=N\cap
L(N_i)$ so $N_i\subseteq L(N_i)$, hence $\dim N_i\leq\dim
L(N_i)$. Since $L$ is a linear operator, $\dim L(N_i)\leq\dim N_i$
holds as well, so $N_i=L(N_i)$.  But $N_i$ is a subspace of $N$, so it
follows that $N_i=\{0\}$. This means that $\dim(N_{i+1}) = \dim(N_i)$ only if $N_i = \{0\}$ so that $\dim(N_{i+1}) \leq \dim(N_i) - 1$ if $\dim(N_i) \neq 0$. This implies our monotonically decreasing upper bound for the dimensions of the $N_i$.

\smallskip

To see that $N_{i+1}\subseteq N_i$: we use induction. We  have
$N_2=N\cap L(N)\subseteq N=N_1$, so $N_2\subseteq N_1$.
If $N_{k+1}\subseteq N_k$, then
$N_{k+2}=N\cap L(N_{k+1})\subseteq N\cap L(N_k)=N_{k+1}$.

\end{proof}

Can one find such a bad example? Yes! Consider the case in which $X =
\Bbb{R}^6$ and

\smallskip

\begin{description}
\item[] \[ L = \left[ \begin{array}{cccccc}
                 0&0&0&0&0&1\\
                 1&0&0&0&0&0\\
                 0&1&0&0&0&0\\
                 0&0&1&0&0&0\\
                 0&0&0&1&0&0\\
                 0&0&0&0&1&0\\
            \end{array} \right] \]

\item[] \[ P = \left[ \begin{array}{cccccc}
                 1&0&0&0&0&0\\
                 0&1&0&0&0&0\\
            \end{array} \right] \]

\end{description}

Then we get that N is given by 

\begin{description}
\item[] \[ N = \left[ \begin{array}{c}
                 0\\
                 0\\
                 x\\
                 x\\
                 x\\
                 x\\
\end{array}
\right]
\]
\end{description}

where the $x$'s can be any value. This gives 


\begin{alignat}{5} N_1 &= \left[ \begin{array}{c}
                 0\\
                 0\\
                 x\\
                 x\\
                 x\\
                 x\\
\end{array}
\right] N_2 &= \left[ \begin{array}{c}
                 0\\
                 0\\
                 0\\
                 x\\
                 x\\
                 x\\
\end{array}
\right] N_3 &= \left[ \begin{array}{c}
                 0\\
                 0\\
                 0\\
                 0\\
                 x\\
                 x\\
\end{array}
\right] N_4 &=  \left[ \begin{array}{c}
                 0\\
                 0\\
                 0\\
                 0\\
                 0\\
                 x\\
\end{array}
\right] N_5 &=  \left[ \begin{array}{c}
                 0\\
                 0\\
                 0\\
                 0\\
                 0\\
                 0\\
\end{array}
\right] \end{alignat}


\section[Optimal Dynamics are generic \textnormal{\\[1ex] \footnotesize  Does a typical problem lead to a sequence of intersections whose dimensions converge to zero as fast as possible? The answer is that in a carefully defined sense YES!}]{Optimal Dynamics are generic}
\label{sec:opt}

We now present and prove a theorem that addresses the point of how
prevalent operators having the optimal reduction property are. In fact
as the title of this section suggests, optimal $L$ are generic. By
generic we mean open and dense (not residual).

\smallskip

In the following  proof, we will use $d_A$ to denote dim($A$). Since
the set of all invertible $\hat{L}$ in $\mathfrak{T}$ is open and of
full measure in $\Bbb{R}^{(T-1)n^2}$ we shall assume that $\hat{L}$ is
invertible throughout the proof. We shall also use the fact that for
invertible $L$,  $ L(A \cap L(B)) = L(A) \cap L^2(B)$.

\bigskip

\begin{thm}[Extended Linear Transverse Intersection Theorem]
\label{thm:eltit}
If the set of operators $\hat{L} = (L_1,L_2,...,L_{T-1})$ is
identified with  $\Bbb{R}^{(T-1)n^2}$ and we define $\mathfrak{T}
\subset \Bbb{R}^{(T-1)n^2}$ to be all those $\hat{L} \in
\Bbb{R}^{(T-1)n^2}$ such that

\[ L_i(N_i) \; \transverse \; N \;\;\; \forall i.\]

Then $\mathfrak{T}$ is open and dense in $R^{(T-1)n^2}$.

\end{thm}

\begin{proof}

\bigskip

We first observe that $\; M_1 \; \transverse \; M_2\Leftrightarrow
d_{P_{M_1^\perp}(M_2)} = \text{min}(d_{M_1^\perp},d_{M_2})$ or
equivalently that $\text{rank}(P_{M_1^{\perp}}\circ P_{M_2}) =
\text{min}(d_{M_1^\perp}$,$d_{M_2})$. A little bit of thought is
enough to convince oneself that $\text{rank}(P_{N^{\perp}} \circ
P_{L(N)}) = \text{rank}(P_{N^{\perp}} \circ L \circ P_N)$. Therefore
we get that $ N \; \transverse \; L(N) \Leftrightarrow
\text{rank}(P_{N^{\perp}} \circ L \circ P_N) = \text{min}(
d_{N^{\perp}}, d_N )$.(Here we have used the invertibility of $L$ to
conclude that $d_N = d_{L(N)}$).  Let us approach the problem a little
more generally. We shall use the fact that $ K \; \transverse \; L(M)
\Leftrightarrow \text{rank}(P_{K^{\perp}} \circ L \circ P_M) =
\text{min}( d_{K^{\perp}}, d_M )$ to show that the set $\mathfrak{T}_*
$ of all $L$ in $\Bbb{R}^{n^2}$ such that $ K \; \transverse \; L(M)$
is open and of full measure. Here $K$ and $M$ are linear subspaces of
$\Bbb{R}^n$.

\bigskip

Define column vectors $p_{\cdot,i}$ for $i = 1,...,d_{K^\perp}$  that
are orthogonal to each other and span $K^{\perp}$. Likewise let
$q_{\cdot,i}$ for $i = 1,...,d_M$ be column vectors orthogonally
spanning $M$. Define the $n$ by $n$ matrices $P$ and $Q$ as follows:

\[   P = \left( \begin{array}{ccccccc}
         p_{1,1} &  p_{1,2} &  \ldots  & p_{1,d_{K^{\perp}}}  & 0 & \ldots  & 0\\ 
         p_{2,1} &  p_{2,2} &  \ldots  & p_{2,d_{K^{\perp}}}  & 0 & \ldots  & 0\\ 
         \vdots  &  \vdots  &  \ddots  & \vdots               &\vdots& &\vdots \\
         p_{n,1} &  p_{n,2} &  \ldots  & p_{n,d_{K^{\perp}}}  &  0 & \ldots & 0
\end{array} \right) 
\label{eq:P}\]

and 

\[   Q      = \left( \begin{array}{ccccccc}
         q_{1,1} &  q_{1,2} &  \ldots  & q_{1,d_M}  & 0 & \ldots  & 0\\ 
         q_{2,1} &  q_{2,2} &  \ldots  & q_{2,d_M}  & 0 & \ldots  & 0\\ 
         \vdots  &  \vdots  &  \ddots  & \vdots     &\vdots& &\vdots \\
         q_{n,1} &  q_{n,2} &  \ldots  & q_{n,d_M}  &  0 & \ldots & 0
\end{array} \right) 
\label{eq:Q} \]

Then

\[   P_{K^{\perp}} = P \circ P^T \]

and 

\[   P_M = Q \circ Q^T, \]

so that we get

\[  P_{K^{\perp}} \circ L \circ P_M = P \circ P^T \circ L \circ Q \circ Q^T. \]

Now we note that 

\[  \text{rank}(P \circ P^T \circ L \circ Q \circ Q^T) = \text{rank}(P^T \circ L \circ Q).\]

\bigskip

To show that $\mathfrak{T}_*$ is open in $\Bbb{R}^{n^2}$, observe that

\[ P^T \circ L \circ Q  = \qquad \begin{array}{cc} 
          & {\tiny \begin{array}{cccc} & d_M & n-d_M & \end{array} } \\ 
            {\tiny \begin{array}{c} d_{K^{\perp}}\\  \\ \\  \\ n -d_{K^{\perp}} \end{array} }
           & \left({\Huge \begin{array}{cc} U_L&0\\0&0 \end{array} } \right) 
             \end{array} \]

and therefore

\[\text{rank}(P^T \circ L \circ Q) = \text{min}( d_{K^{\perp}}, d_M ) \Leftrightarrow
   U_L \; \text{is full rank}. \]

Note that $U_L$ is a continuous function of $U$, i.e.
$U_L:\Bbb{R}^{n^2} \rightarrow \Bbb{R}^{d_M \cdot d_{K^{\perp}}}$  is
continuous (actually smooth). Assume without the loss of generality, that $d_{K^{\perp}} \geq  d_M$.   Let $\phi_{d_M}^{d_{K^{\perp}}}(U_L)$
be the $d_M$ dimensional measure of regions in $\Bbb{R}^{K^{\perp}}$
applied to the parallelepiped with edges equal to the columns of
$U_L$. Then $\mathfrak{T}_*$ is precisely equal to
$(U_L)^{-1}((\phi_{d_M}^{d_{K^{\perp}}})^{-1}(\Bbb{R}\setminus\{0\}))$.
Since both $U_L$ and $\phi_{d_M}^{d_{K^{\perp}}}$ are continuous, we
have that $\mathfrak{T}_*$ is open.

\bigskip

To show that $\Bbb{R}^{n^2} \setminus \mathfrak{T}_*$ has zero
$n^2$-dimensional Lebesgue measure, we first introduce a change of
coordinates. Define $\hat{P}$ to be an orthogonal matrix obtained
by filling in the zero columns of $P$ appropriately. Obtain $\hat{Q}$
from $Q$ analogously. Then

\begin{eqnarray*}
P^T\circ L\circ Q =& P^T\circ\hat{P}\circ\hat{P}^T\circ L\circ\hat{Q}\circ\hat{Q}^T \circ Q \\
 =& \left(\begin{array}{cc} I_{d_{K^{\perp}}} & 0 \\ 0 & 0 \end{array} \right) \circ
 \hat{L}_L \circ \left( \begin{array}{cc} I_{d_M} & 0 \\ 0 & 0 \end{array} \right)\\
 =& \hat{L}_L(\text{ul})\\
 =& \text{ upper left} \; (d_{K^{\perp}} \times d_M)\; \text{block of} \; \hat{L}_L
\end{eqnarray*}

where $I_{\zeta}$ is the identity matrix of dimension $\zeta$
and we have set $\hat{L}_L = \hat{P}^T\circ L\circ\hat{Q}$. So 

\[\begin{split}
\mathfrak{T}_* =& \{ L | \; \text{rank}( \hat{L}_L(\text{ul}))  = \text{min}(d_{K^{\perp}},d_M)\; \}\\
      =&\{ L | \;  \hat{L}_L(\text{ul}) \; \text{is full rank} \}
\end{split} 
\]

Since $\hat{P}$ and $\hat{Q}$ are orthogonal, we have that the
$n^2$ dimensional Lebesgue measure of

\[ \mathfrak{T}_* = \{ L | \;  \hat{L}_L(\text{ul}) \; \text{is full rank} \}
\]

and

\[ \hat{\mathfrak{T}}_* \equiv \{ \hat{L} | \; \hat{L}(\text{ul}) \; \text{is full rank} \}
\]

are equal.

\bigskip

We now prove that $\Bbb{R}^{n^2} \setminus \hat{\mathfrak{T}}_*$ has measure zero. Any 
 $\hat{L}$ can be written as the block matrix with dimension of $\hat{L}_{ul}$ being $d_{K^{\perp}} \times d_M$.

\[ \left(  \begin{array}{cc}  \hat{L}_{ul} & \hat{L}_{ur} \\
                              \hat{L}_{ll} & \hat{L}_{lr} \end{array} 
   \right) \]

We can write $\hat{L}$ out in terms of elements as

\[ \hat{L} = \left( \begin{array}{cccc} 
                       \hat{l}_{11} & \hat{l}_{12} & \ldots & \hat{l}_{1n}\\
                       \hat{l}_{21} & \hat{l}_{22} & \ldots & \hat{l}_{1n}\\
                        \vdots        &    \vdots      & \ddots &   \vdots      \\
                       \hat{l}_{n1} &  \hat{l}_{n2} & \ldots & \hat{l}_{nn}
                       \end{array}
                \right).
\]
 
Now assume that $d_{K^\perp} \geq d_M$. Identify  $\Bbb{R}^{n^2 -
d_{K^\perp}\cdot d_M + (d_M - 1) \cdot ( d_{K^\perp} + 1)}$   with the
elements of $f_{ij}$ of an $n \times n$ matrix with the elements
$f_{d_M,d_M}\; ,\;f_{d_M+1,d_M} \; , \; ... \; , \;
f_{d_{K^\perp},d_M}$ removed.

Next, define a mapping of 
$\phi : \Bbb{R}^{n^2 - d_{K^\perp}\cdot d_M + (d_M - 1) \cdot ( d_{K^\perp} + 1)} \rightarrow
        \Bbb{R}^{n^2}$ by  

\[
  \hat{l}_{ij} = \begin{cases} f_{ij} & \text{for $(i,j) \; \not\in$ \{$i \leq d_{K^\perp}$ and $j = d_M$\} } \\
                  \sum f_{kd_M}f_{ik} & \text{for $(i,j) \; \in$ \{$i \leq d_{K^\perp}$ and $j = d_M$\} }  
                   \end{cases}
\]

which has, as its image precisely those matrices in which
$\hat{L}_{ul}$ has column $d_M$ that is the linear combination of
the first $d_M - 1$ columns. If we redefine our mapping to get a
series of  completely analogous mappings, each of which has a column
of $\hat{L}_{ul}$ being dependent on the other columns
of $\hat{L}_{ul}$, then we end up with $d_M$ such maps. Since each of
these maps are smooth, and singular (the rank of the derivative is not
equal to the dimension of the image space) Sard's theorem ~\cite{milnor-TDV} tells us
that the $n^2$-dimensional measure of the image of each map is
zero. Therefore the union of the images also has measure zero. But
this union is exactly $\Bbb{R}^{n^2} \setminus
\hat{\mathfrak{T}}_*$. Since the case of  $d_{K^\perp} < d_M$ is
completely analogous, we have now shown that $\mathfrak{T}_*$ is open
and of full measure in $\Bbb{R}^{n^2}$.

\bigskip

To complete the proof we let the operator
change at each step so that $x_2 = L_1(x_1)$ , $x_3 =
L_2(x_2)$ , and so on. Now we have an extended operator $\tilde{L} =
(L_1,L_2,...,L_{T-1}) \in (\Bbb{R}^{n^2})^{T-1}$. We will show that the set
$\mathfrak{S} \equiv \{ \tilde{L} | \; N \transverse L_i(N_i) \;
\text{for}\; i = 1,2,...,T-1\}$ is open and dense in
$(\Bbb{R}^{n^2})^{T-1}$. 

\bigskip

Define $\Bbb{C}_1$ to be the open subset of full measure in
$\Bbb{R}^{n^2}$ whose members, $L_1$, satisfy $N \; \transverse \;
L_1(N_1)$.  Choose a countable subset, $\Bbb{D}_1$, which is dense in
$\Bbb{C}_1$. Now for each element $D_1^k$ of $\Bbb{D}_1$, define
$\Bbb{C}_2^k$ to be the open full measure subset of $\Bbb{R}^{n^2}$
such that $L \in \Bbb{C}_2^k$ implies that $N \; \transverse \; L(N
\cap D_1^k(N_1))$, or equivalently, $N \; \transverse \; ( L(N) \cap
L\cdot D_1^k(N_1) )$. Define  $\Bbb{C}_2 \equiv \bigcap_k
\Bbb{C}_2^k$. Since $\Bbb{C}_2$ is dense in $\Bbb{R}^{n^2}$ we can
pick a countable $\Bbb{D}_2 \subset \Bbb{C}_2 $ that is also dense in
$\Bbb{R}^{n^2}$. We have that $L_1 \in \Bbb{D}_1 $ and $L_2
\in\Bbb{D}_2 $ implies that $N \; \transverse \; L_1(N_1)$ and $N \;
\transverse \; L_2(N \cap L_1(N_1))$. Continuing this process we
obtain $\Bbb{D}_i$ for $i=1,2,...,T-1$ such that
$(L_1,L_2,...,L_{T-1}) \in  \Bbb{D}_1 \times \Bbb{D}_2 \times
... \times  \Bbb{D}_{T-1}$ implies that all the intersections are
transverse, i.e. that $N \; \transverse \; L_i(N_i)$ for
$i=1,2,...,T-1$. We have therefore found a subset of $\mathfrak{S}$
which is dense in $(\Bbb{R}^{n^2})^{T-1}$. 

\bigskip

Now we show that $\mathfrak{S}$ is open. The requirement that each of
the intersections are transverse is equivalent to the requirement that

\[\begin{split}
\text{dim}(N \cap L_1(N)) &= d_N + d_N - n \\
\text{dim}(N \cap L_2(N) \cap L_2L_1(N)) &= 3d_N  - 2n \\
\text{dim}(N \cap L_3(N) \cap L_3L_2(N) \cap L_3L_2L_1(N)) &= 4d_N - 3n \\
      \vdots         \qquad   &=   \qquad      \vdots               \\
\text{dim}(N \cap L_{T-1}(N) \cap L_{T-1}L_{T-2}(N) \cap ... \cap L_{T-1}...L_1(N)) &= Td_N-(T-1)n. 
\end{split}\]

which in turn is equivalent to a requirement involving orthogonal
complements, specifically that

\[  {\tiny 2(n-d_N)} \overset{\large n}{ \left[\begin{array}{c} rN^\perp \\ rN^\perp \circ L_1^{-1}\end{array} \right]} \]

\[  {\tiny 3(n-d_N)} \overset{\large n}{ \left[\begin{array}{c} rN^\perp \\ rN^\perp \circ L_2^{-1} \\  rN^\perp \circ L_1^{-1} \circ L_2^{-1} \end{array} \right]} \]

\[ \qquad \vdots \]

\[  {\tiny T(n-d_N)} \overset{\large n}{ \left[\begin{array}{c} rN^\perp \\ rN^\perp \circ L_{T-1}^{-1} \\ \vdots \\ rN^\perp \circ  L_{1}^{-1} \circ ... \circ  L_{T-2}^{-1} \circ L_{T-1}^{-1}  \end{array} \right]} \]

all have full rank where $rN^\perp$ is the matrix with rows equal to
independent n-dimensional vectors spanning the linear subspace
$N^\perp$. This last set of expressions follows from the fact that if
$M$ and $K$ are linear subspaces of $\Bbb{R}^n$ then $(M \cap K)^\perp
= \text{span}(M^\perp , K^\perp)$ so that the matrices immediately above
have rank $2(n-d_N),3(n-d_N),...,T(n-d_N)$ iff the the previous
intersections have dimensions $2d_N-n, 3d_N-2n,..., Td_N-(T-1)n$
respectively. But this last expression can be seen to be exactly
those matrices which satisfy the equations

\[ \phi_{2(n-d_N)}^n(\text{rows from first matrix}) \neq 0 \]
 
\[ \phi_{3(n-d_N)}^n(\text{rows from second matrix}) \neq 0 \]

\[   \qquad \vdots                \]

\[ \phi_{T(n-d_N)}^n(\text{rows from T - 1st matrix}) \neq 0 \]

\bigskip

where $\phi_j^l(vectors)$ measures the $j$-dimensional volume of the
the parallelepiped spanned by the $vectors$ in $\Bbb{R}^l$. But since the
inverse operation is continuous on the set of invertible matrices and
these volume functions are smooth, we have that the set of
$(L_1,L_2,...,L_{T-1})$ having the full rank property is open. Thus,
$\mathfrak{S}$ is open in $(\Bbb{R}^{n^2})^{T-1}$.

\end{proof}

\bigskip

In the next section we conjecture an approach to the nonlinear case
which uses the above derivation, but before we do this we show that
the linear case can actually be made significantly simpler.

\bigskip

\begin{thm}[Linear Transverse Intersection Theorem]

If the set of operators $L$ is
identified with  $\Bbb{R}^{n^2}$ and we define $\mathfrak{T}
\subset \Bbb{R}^{n^2}$ to be all those $\L \in
\Bbb{R}^{n^2}$ such that

\[ L(N_i) \; \transverse \; N \;\;\; \forall i.\]

Then $\mathfrak{T}$ is open and dense in $R^{n^2}$.

\end{thm}

\begin{proof}

As was seen in the proof of the previous theorem, the transversality 
requirement and the fact that we are considering the case of $L_i = L$ for all $i$, reduces to

\[  {\tiny 2(n-d_N)} \overset{\large n}{ \left[\begin{array}{c} rN^\perp \\ rN^\perp \circ L^{-1}\end{array} \right]} \]

\[  {\tiny 3(n-d_N)} \overset{\large n}{ \left[\begin{array}{c} rN^\perp \\ rN^\perp \circ L^{-1} \\  rN^\perp \circ L^{-2} \end{array} \right]} \]

\[ \qquad \vdots \]

\[  {\tiny T(n-d_N)} \overset{\large n}{ \left[\begin{array}{c} rN^\perp \\ rN^\perp \circ L^{-1} \\ \vdots \\ rN^\perp \circ  L^{-(T-1)} \end{array} \right]} \]

all having full rank. But this is equivalent to another full rank condition
as follows. Let $cN^\perp$ be the transpose of $rN^\perp$. In other words,
while $rN^\perp$ are the orthogonal row vectors that span compliment of the 
null space, $cN^\perp$ are the column vectors that span the ``same'' space.
This condition is equivalent to the following. For a dense and open set of
 $L$:

\[ \text{dim}(S_k \equiv \text{span}(cN^\perp,L\circ cN^\perp,...
    ,L^{k-1}\circ cN^\perp)) = \text{min}(k\cdot d_{cN^\perp},n) \]

That is, for an open and dense set of $L$ the sequence of subspaces
$S_k$ generated by the iterates  $L^j\circ cN^\perp$ $j = 1,2,...,k-1$, are
of maximal dimension.

\bigskip

The proof of this fact is well known. For lack of a reference I give a 
proof here. Without loss of 
generality let $cN^\perp$ be the k left most columns of the
$n\times n$ identity matrix. Then the dimension of $S_k$ is the 
rank of the matrix formed by taking the k left columns of I,
followed by the k left columns of L, followed by the k left
columns of $L^2$, and so on.

\bigskip

Pick the upper left matrix minor and compute it's determinate this
will give a polynomial in $l_{11},l_{12}, l_{21},... l_{nn}$ which we
want to show is nonzero except on an open and dense set. As long as
the polynomial is nonzero at one point then we are done since this
implies that the set of zeros occupies a submanifold of $R^{n^2}$ that
is at most, $n^2 - 1$ dimensional. (So we even have more ... the set
of ``good'' $L$ has full measure and is open!)

\bigskip

To show that the determinant of the upper left matrix minor is nonzero
at a point, we consider $L$ = permutation matrix that shifts
everything to the left k clicks. This gives us the identity matrix in
the upper left matrix minor and so the determinant in question
evaluates to 1!

\end{proof}

{\bf Remark:} As the above proof shows, $\mathfrak{T}$ is in fact open
and full measure. This improves the result from one stated which
implies stable approximation by $\hat{L}$ having the optimal reduction
property, to one that implies this AND the improbability of
non-optimal reduction. Further improvements would involve the
characterization of $\mathfrak{T}_\epsilon \equiv \mathfrak{T} \cap \{
L \in \Bbb{R}^{n^2}| \text{cond}(L) \;\; < \;\; 1/ \epsilon \}$.

\bigskip

{\bf Remark:} The above proof also works with slight modifications to
prove theorem \ref{thm:eltit}, but the proof given there leads to a
conjectured proof for the case of nonlinear dynamics, and so seems
more useful even if it is more cumbersome.

\bigskip

\section[Extension to the Nonlinear Case \textnormal{\\[1ex] \footnotesize So far everything we have developed has been for the case in which the underlying dynamics are linear. What can we say about the nonlinear case? We use our hands (waving them about) to say something which we make respectable by  formulating two conjectures and sketching how we think it should go.}]{Extension to the Nonlinear Case}

\label{sec:nonlin}

The above theorem is extendable to the nonlinear
case as follows. Actually, we are conjecturing such an extension in what
follows.  It should be noted that the nonlinear ``extension'' does not imply the  linear theorem. 

\bigskip

In this section the state space (object space) will be $M$, a compact
manifold of dimension $m$. The projection operator will be a smooth
function $P:M \rightarrow \Bbb{R}^d$ and the dynamics will be given by
$F_i$'s which map $M$ to itself diffeomorphically. Instead
of using linearity to get a fixed null space, we find the the ``null''
space we are now interested in is a level set of $P$. These sets can 
change in nontrivial ways as the point in the range of $P$ changes. We
now want to know about intersections of these ``null'' sets with 
images of other of the ``null'' sets under $F$.

\bigskip

The transversality theorem found on page 74 of ~\cite{hirsch-DT}
implies that the set $\mathfrak{T} $ of $f \in C^r(M)$ which map
a submanifold $K$ of $M$ back into $M$ to intersect another submanifold $N$
transversely, is open and dense. (For compact $M$ the topology is
nice, see ~\cite{hirsch-DT} for details.) What we need is a bit more
complicated.

\medskip

Let  $\mathfrak{F}$ be the quasi-stratification of $M$ into the
level-sets $\mathfrak{F}_{x}$, $x \in R^d$ of $P$. Let N be the union
of a finite number of (not-necessarily injectively) immersed 
submanifolds of dimension $ \leq n$ whose self-intersections are transverse in
the sense that the tangent spaces of the ``participants'' in the
intersection span the largest possible subspace of $T_{i_x}M$ where
$i_x$ is a point of self-intersection.

\medskip

\begin{con} 

For an open and dense set of $F$ $\in  D^{\infty}(M) $

\[  dim(I_x) \leq \text{max}(n-d,0) \;\; \forall x \in \Bbb{R}^d \]

where $I_x \equiv (F(N) \cap \mathfrak{F}_{x})$ and,  

\[  I_x  \;\; \mathrm{is\ a\ finite\ union\ of\ stably\ 
immersed\ submanifolds\ }  \].

\end{con}

\medskip

This  permits us to conclude that, for any initial point $x_0 \in X$
and dense $\tilde{F} = (F_1,...,F_{T-1})$, the intersection obtained by
$T = \lceil m/d \rceil$ measurements will have dwindled to a finite
set of points, call it $S$. The next measurement ($T =\lceil m/d \rceil+1$) 
 will generically have
precisely one point in the intersection of the set $S$ and the level
set corresponding to the next measurement. (We use the same argument
as we used in the linear case to get a product of dense sets
$\mathfrak{D}_1 \times ... \times \mathfrak{D}_{T-1}$.)

\medskip

Now, if indeed the $\lceil m/d \rceil+1$``th'' intersection
 is a single point, then the
mapping $G:x \in M \rightarrow (Px,PF_{1}x,PF_{2}F_{1}x, ... ,
PF_{T-1}...F_{1}x)$ has only one point in the inverse image of
$G(x_0)$. Since the point we have ``found'' (the conjecture above)
comes from stable intersections this should guarantee that the point
$G(x_0)$ in fact has a neighborhood in which $G$ is invertible. This
should in turn guarantee that there is an open neighborhood $B_{\epsilon}$
of  $G$ in $C^{\infty}(M)$ such that $H \in B_{\epsilon}$ implies
$H^{\leftarrow}(H(x_0))$ is a single point. We then use the fact that a small
neighborhood of $\tilde{F}$ maps into this small neighborhood under the
mapping $\tilde{J} \rightarrow (P,PJ_1, ... ,PJ_{T-1}...J_{1})$  for $J \in D^{\infty}(M)$.

\bigskip

We have arrived at our second conjecture.

\begin{con}

If $M$ is a compact smooth manifold of dimension $m$, $P$ is a  smooth
function mapping $M$ to $\Bbb{R}^d$, $x_0$ is a particular point in
$M$, $T = \lceil m/d  \rceil + 1$, $\mathfrak{D} \equiv
(\Bbb{D})^{T-1}$ is the $T-1$-fold product of the space of smooth
diffeomorphisms from $M$ to $M$ ($D^{\infty}(M)$) and we define
$H_\delta \equiv  (P,P\delta_1,...,P\delta_{T-1}...\delta_1): M \rightarrow
\Bbb{R}^{dT}$, where $\delta= (\delta_1,...,\delta_{T-1})
\in \mathfrak{D}$ , then the set $\mathfrak{O}$ of $\delta \in 
\mathfrak{D}$ such that

\[ H_{\delta}^{\leftarrow}(H_{\delta}(x_0)) = \{x_0\}  \]

is open and dense in $\mathfrak{D}$.
   
\end{con}

\medskip

\section[Relation to Known Results \textnormal{\\[1ex] \footnotesize Of course even though this work was new to us (and therefore fun!) there is nothing new under the sun and we find other work that is very similar. Our hope is that the differences will be useful to the readers and that the readership will be much wider than the few experts who would find the contents not surprising.}]{Relation to Known Results}
\label{sec:RKR}

The results obtained above are known as observability results in
control theory and phase space reconstructions (delay coordinate
embedings) in dynamical systems. Our results are different in that we
consider variations of the dynamics with the observation function kept
fixed  whereas other results either assume that the dynamics are fixed
and the observation function changes or that both the dynamics and the
observation function is variable, see
~\cite{takens-1981,sauer-E,aeyels-1981,balde-1998,stark-1999}. While
Aeyels ~\cite{aeyels-1981} does consider the case where the
observation function is fixed and the dynamics are variable he does so
for vector fields (not maps). He is also looking at the case where he
wants all initial points to be recoverable from the sequence of
measurements  and there requires 2n+1 , 1-dim measurements to recover
the n-dim initial points. Similar comments apply to the comparison to
Stark's more recent paper ~\cite{stark-1999}. Our minimal reduction
theorem is a more precise version of the well known theorem in control
theory that states that if the observability matrix is not full rank
then no number of measurements can give you full information on the
state of the system and if it is full rank then you need at most n
measurements (of any dimension) of a system that has an n-dimensional
state space.

\section[Numerical Examples \textnormal{\\[1ex] \footnotesize Here we demonstrate the observability of a particular linear system as well as some of the problems associated with these type of efforts.}]{Numerical Examples}
\label{sec:num}

We now give two examples in which we apply the technique described
in section \ref{sec:lin} to invert simulated sequences of
(noiseless) radiographs using one view.  Our purpose in this section is
simply to demonstrate the procedure.  We assume our object lies within
a $10\times 10$ pixelation and has constant density within each pixel,
so the object space $X$ is $\Bbb{R}^{100}$.  We use the same initial
condition with two different linear operators $L_1$ and $L_2$ which we
describe below.  In each case, the projection $P$ sums the values down
the columns of the pixelation.  We use, in a sense, the largest
parameterization of our object space; namely, we assume nothing about
the object and seek to determine the value in each pixel.  At the end
of this section, we comment briefly on the poor numerical conditioning
of these problems and indicate first steps taken to improve the
numerics.  This is a subject of current study.

To reiterate the procedure, first choose a basis $\{b_i\}_{i=1}^{100}$
for $X$.  With $L$ and $P$ representing the dynamics and projection
operators respectively, we build a $100\times 10t$ matrix $E$ where the
$i$th row of $E$ is
$(Pb_i,PLb_i,PL^2b_i\ldots,PL^{t-1}b_i)=\tilde{P}_t\tilde{b}_i$.
As soon as $t$ is large enough so that $ \Rank E=100$, we have a unique
solution $x$ for the equation $xE=\tilde{d}^*$.  Since we know the dynamics
$L$, we can then reconstruct the sequence
$\tilde{x}^*=(x,Lx,L^2x,\ldots,L^{t-1}x)$.

In each example, we chose the canonical basis $\{e_i\}_{i=1}^{100}$
for $X$ where $e_i(j)=\delta_{ij}$ $(1\leq j\leq 100)$.
The first linear operator $L_1$ can be described as a combination
of a diffusion and a shift. The effect of $L_1$ is pictured below for
various times $t$.
\begin{figure}[h]\label{fig:fig1}
\centering
\mbox{
\subfigure[$x$]
{\epsfig{file=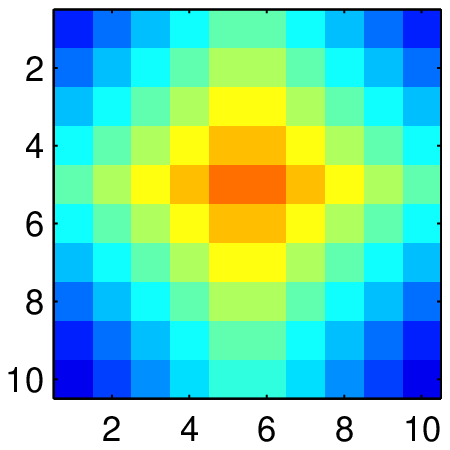, height=1.05in, width=1.05in}}
\hspace{0.3cm}
\subfigure[$L_1^2x$]
{\epsfig{file=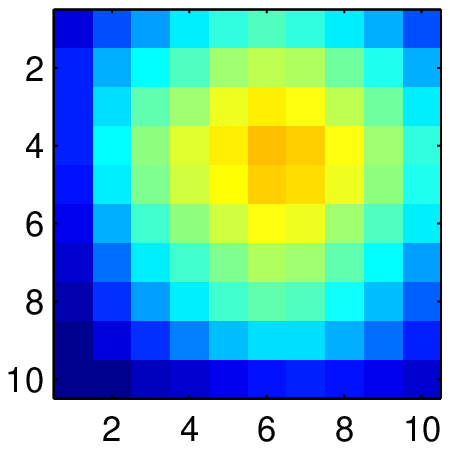, height=1.05in, width=1.05in}}
\hspace{0.3cm}
\subfigure[$L_1^5x$]
{\epsfig{file=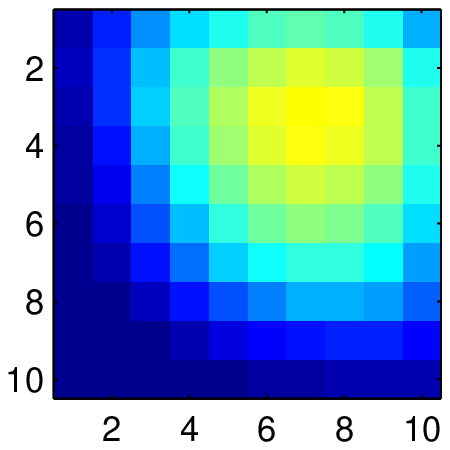, height=1.05in, width=1.05in}}
\hspace{0.3cm}
\subfigure[$L_1^9x$]
{\epsfig{file=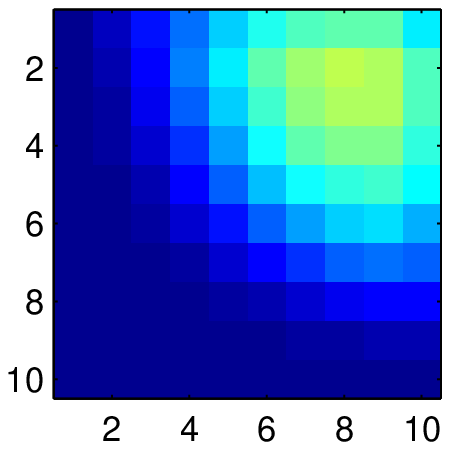, height=1.05in, width=1.05in}}}
\caption{Initial condition and $L_1^tx$ for $t=2,5,9$.}
\end{figure}

\noindent

Here, the rank of $E$ increased by 10 each time step, so we achieved
$ \Rank E=100$ in the minimal number of steps and were able to solve for
the initial condition $x$.

Our second operator, $L_2$, was a diffusion operator where the diffusion
coefficient varied over the pixelation. Namely, the diffusion coefficient
in the $i,j$-pixel was $\frac{1}{5}(i^3j^210^{-5})^{1/4}$ (so the rate of
diffusion was greatest in the lower right corner of the pixelation and
was least in the upper left corner).  The effect of $L_2$ is pictured below
for a few times $t$.

\begin{figure}[h]\label{fig:fig2}
\centering
\mbox{
\subfigure[$x$]
{\epsfig{file=images/init_cond.ps, height=1.05in, width=1.05in}}
\hspace{0.3cm}
\subfigure[$L_2^2x$]
{\epsfig{file=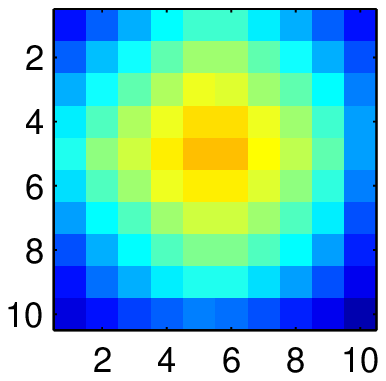, height=1.05in, width=1.05in}}
\hspace{0.3cm}
\subfigure[$L_2^5x$]
{\epsfig{file=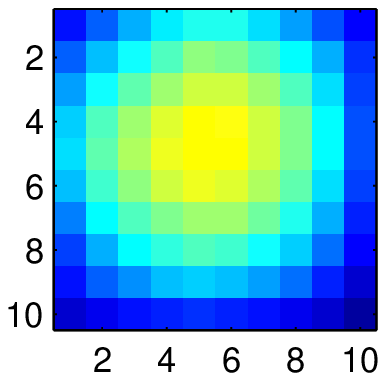, height=1.05in, width=1.05in}}
\hspace{0.3cm}
\subfigure[$L_2^9x$]
{\epsfig{file=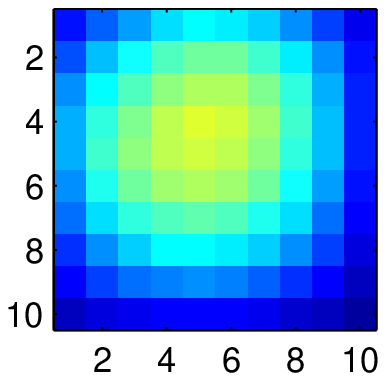, height=1.05in, width=1.05in}}}
\caption{Initial condition and $L_2^tx$ for $t=2,5,9$.}
\end{figure}

\noindent
In this example, the rank of $E$ again increased by 10 each step reaching
100 after 10 steps.  Pictured below are the initial condition $x$ and the
reconstructions obtained by using the data sequence
$(Px^*,PL_2x^*,PL_2^2x^*,\ldots,PL_2^tx^*)$ for $t=9,14$.
\begin{figure}[h]\label{fig:fig3}
\centering
\mbox{
\subfigure[$x^*$]
{\epsfig{file=images/init_cond.ps, height=1.05in, width=1.05in}}
\hspace{0.5in}
\subfigure[$x, t=9$]
{\epsfig{file=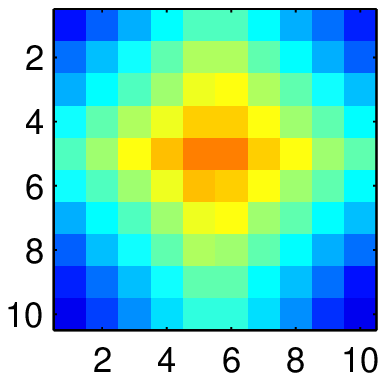, height=1.05in, width=1.05in}}
\hspace{0.5in}
\subfigure[$x, t=14$]
{\epsfig{file=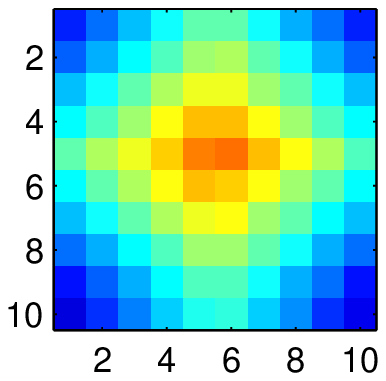, height=1.05in, width=1.05in}}}
\caption{Initial condition and reconstructions using $L_2$ for $t=9,14$.}
\end{figure}

With regard to the numerical conditioning of these problems, we note
that the condition numbers of the matrices $E$ constructed using $L_1$
and $L_2$ were on the order of $10^{12}$ and $10^{11}$ respectively.
By running the dynamics longer than the number of time steps required
to achieve full rank, we were able to reduce the condition number in
both cases.  Namely, using $L_1$ for 15 time steps reduced the
condition number of $E$ to $10^{11}$.  But with $L_2$, using 12
time steps reduced the condition number of $E$ to $10^9$, and at 15
time steps the condition number reduced to $10^8$.  In both cases,
extending beyond 15 steps gave no significant improvement.

\section[Summary and Discussion\textnormal{\\[1ex] \footnotesize This is the part (along with the abstract!) that will be read by most readers -- given how busy most researchers are and how they tend to read papers. Therefore we have implanted devious traps that insure those reading this section will inevitably get sucked into reading the whole paper ... (just joking!).}]{Summary and Discussion}
\label{sec:sum}

In this paper we examined the use of dynamics in the inversion
of projection data obtained at a sequence of times. The main results
confirm that for any fixed measurement projection and generic
dynamics, we can simply combine the number of measurements into one
large super-measurement which we invert to obtain the state we are
trying to reconstruct. A following paper will deal with some
aspects of the stochastic or noisy case of reconstruction from
projections using dynamics.

\bigskip

What we have established is only a first step in the direction
leading to the fruitful combination of dynamics and measured data.
Many variants of the proposed underlying tomography problem lead  to
the same abstract problem that we have begun to examine. For  example,
if the dynamical propagator f is known up to a set of parameters, the
resulting abstract inverse problem is identical to the one stated
above. There is still a state space one is trying to observe, only now
it is $n+p$ dimensional where $p$ is the dimension of the  parameter
space. The projection (measurement) function now defines $n+p-d$
dimensional level sets that are mapped forward by the dynamics as
before. We end up needing $ \lceil (n+p)/d  \rceil$ measurements. It
seems to us that there are at least two important directions to go
next. One is the examination of the present formulation in the
presence of noise. This will bring us much closer to ``real''
situations in that the prevalence of noise makes certain problems
which are well posed in the noiseless case, ill-posed in the presence
of noise. The second direction is the attack of a very
carefully chosen concrete problem involving dynamics that we
understand analytically or at least numerically. This will invariably
involve certain toy-like characteristics which should nevertheless
be useful for the approach to the large, more realistic problems.

\bigskip

Natural questions that arise include:

\bigskip

\begin{itemize}
\item How is the problem of reconstruction from a sequence of
projections related to the reconstruction of a 3-dim object from a
spatial sequence of slices? This arises when one wants to
interpolate a set of CAT scans to generate a 3-dim density image.  A
related problem arises when one is trying to compress a video movie
by using some clever interpolation in the uncompress process.
\item How do we do (algorithmically, efficiently) the reconstruction
in the case of nonlinear  dynamics. Actually, even the linear case,
while conceptually simple,  is not easy in the case of high dimension
and ``noiseless'' (only computationally induced uncertainties and
approximation induced uncertainties). The difficulty is that
extraordinary large condition numbers are pervasive and so any error ,
like roundoff, soon overwhelms you. We will begin to address these
issues in the next paper which looks at reconstruction using dynamics
in the presence of uncertainty.
\item If one has a set of measurements, how does one use the
knowledge of the underlying dynamics (incomplete maybe) and the
freedom to choose the object (state space) parameterization in such a
way as to get a well posed inverse problem with as little as possible
``wasted'' measured information.  That is, How does one use all the
prior information and measured information in the generation of the
final reconstruction. Even if we are using all our information, there
are different ways of distributing remaining uncertainty about the
reconstructed object. At each for a given data set what sort of
different parameterizations give rise to this no wasted information
situation?
\item In the high dimensional case, the questions asked in the 
previous bullet are incredibly hard to answer. What approximate
answers can be generated? Can we tell how far from the optimum we
are? for example, can we obtain bounds on the amount of information
that our parameterization/reconstruction/use-of-priors wastes?
\item Suppose we do the whole analysis with $\epsilon$ fattened null
spaces. In this case, what sort of volume do we get for the final
intersection (which before was just one point)? This is along the same
lines as the first remark at the end of section ~\ref{sec:opt}.
\end{itemize}

\section[Acknowledgments\textnormal{\\[1ex] \footnotesize Thanks to
 all those with whom we have had helpful discussions and all those who
 have forked over the bucks that made the research
 possible.}]{Acknowledgments}
\label{sec:ack}

This research was supported by the University of California under
contract to the Department of Energy, Contract \# W-7405-ENG-36. We
want to thank Tim Sauer for enlightening comments. The work was made
possible by the thought conducive environment of our team devoted to
difficult inverse problems. The members of that team are Dominic
Cagliostro (team leader), Kevin Buescher, James Howse, John Pearson,
Ed Mackerrow, and ourselves. We are thankful that John Pearson made an
``irritating'' remark that turned out to be inspiring. Finally, one of us
(Kevin R. Vixie) thanks Andrew M. Fraser for inspiration and patience.

\bibliographystyle{plain}
\bibliography{/home/vixie/projects/the_bib}

\end{document}